\documentclass[12pt]{amsart}
\usepackage{latexsym}
\numberwithin{equation}{section}
\newtheorem{theorem}{Theorem}[section]
\newtheorem{lemma}[theorem]{Lemma}

\newtheorem{corollary}[theorem]{Corollary}

\theoremstyle{definition}
\newtheorem{definition}[theorem]{Definition}
\theoremstyle{remark}

\newtheorem{example}[theorem]{Example}

\newcommand\codim{\operatorname{codim}}
\newcommand\Ann{\operatorname{Ann}}
\newcommand\Supp{\operatorname{Supp}}
\newcommand\Proj{\operatorname{Proj}}
\newcommand\Ass{\operatorname{Ass}}
\newcommand\reg{\operatorname{reg}}
\newcommand\Tor{\operatorname{Tor}}
\newcommand\Hom{\operatorname{Hom}}
\newcommand\Ext{\operatorname{Ext}}

\newcommand\depth{\operatorname{depth}}
\newcommand{\xx}{\underline x}
\newcommand{\yy}{\underline y}
\newcommand{\zz}{\underline z}
\newcommand{\HH}{H_{\mathfrak m}}
\newcommand{\lu}{\underline l}
\newcommand{\lf}{\underline f}
\newcommand{\s}{\; | \;}

\begin{document}

\author[P.~Schenzel]{Peter Schenzel}
\title[Koszul homology and syzygies]
{Applications of Koszul homology to numbers of generators
and syzygies}

\address{Martin-Luther-Universit\"at Halle-Wittenberg,
Fachbereich Mathematik und Informatik, D --- 06 099 Halle, Germany}

\email{schenzel@mathematik.uni-halle.de}


\begin{abstract}
Several spectral sequence techniques are used
in order to derive information about the structure
of finite free resolutions of graded modules.
These results cover estimates of the minimal number
of generators of defining ideals of projective
varieties. In fact there are generalizations of a
classical result of Dubreil. On the other hand there
are investigations about the shifts and the dimension
of Betti numbers. To this end there is a local
analogue of Green's considerations developed in
\cite{mG}.
\end{abstract}
\subjclass{Primary 14 B 15, 14 M 05; Secondary 13 D 45, 14 M 07}

\maketitle


\section{Introduction}

Let $X \subset \mathbb P^n_K$ be an algebraic variety,
$K$ an algebraically closed field. Let $\mathcal J_X$
denote the ideal sheaf of $X.$ Then $\mathcal J_X$ admits
a finite minimal resolution \begin{displaymath}\mathcal F_{\bullet} : 0 \to
\mathcal F_s \to \ldots \to
\mathcal F_i \to \ldots \to \mathcal F_1 \to \mathcal J_X \to 0,
\end{displaymath}where $\mathcal F_i \simeq \oplus_{j\in \mathbb Z}
\mathcal O^{i_j}(-j).$ Here only finitely many $i_j$ are non-zero.

The resolution $\mathcal F_{\bullet}$ reflects several
geometric and arithmetic properties of $X.$ For
instance, the length $s$ of $\mathcal F_{\bullet}$ satisfies
$s \geq \codim X.$ The equality characterizes when
$X$ is arithmetically Cohen-Macaulay. On the other hand
$$
\reg M :=
\max \{ j - i \s j \in \mathbb Z \mbox{ and } i_j \neq 0 \}
$$
is called the Castelnuovo-Mumford regularity of $X.$
It is determined by the vanishing of the cohomology $H^i(X, \mathcal J_X(n)),$
see e.g., \cite{dM}. The rank of $\mathcal F_1$ is determined by
the minimal number of generators $\mu (I_X)$ of $I_X,$
the saturated ideal of $X$ in $R = K[x_0,\ldots,x_n].$

There is a classical result by P.~Dubreil, see
\cite{pD}, that for $X$ a set of points in the projective
plane $\mathbb P^2_K$ it follows that
$$
\mu (I_X) \leq a(I_X) + 1,
$$
where $a(I_X)$ denotes
the smallest degree of a hypersurface that contains $X.$
There is an extension of this result to the case of $X$
arithmetically Cohen-Macaulay of codimension two, see
\cite{DGM} or \cite{pS1}. More recently H.~Martin and
J.~Migliore extended Dubreil's Theorem to $X$ a locally
Cohen-Macaulay scheme, see \cite[Theorem 2.5]{MM}. One
of the main points of the present paper is an extension
of their result to an arbitrary scheme $X \subset \mathbb P^n_K.$
In fact it turns out that
$$
\mu (I_X) \leq a(I_X) + 1 + C(I_X)
$$
for
a certain correcting term $C(I_X),$ see \ref{4.1} and \ref{4.3}.

The integer $C(I_X)$ is determined by the dimensions of
the Koszul homology modules of $H^i_{\ast}(X, \mathcal O_X)$ with
respect to a certain system of linear forms. These kinds
of invariants have been considered by M.~Green in his fundamental paper
\cite{mG}. In fact
we develop a complete local analogue of these modules
for any finitely generated graded $R$-module $M.$ These
invariants
$$
H_i(\lu ; \HH^j(M)),\quad i, j \in \mathbb Z,
$$
$\lu = l_1, \ldots, l_r, $
a system of generic linear forms,
are graded $R$-modules of finite length. We call them
Green modules of $M$ with respect to $\lu.$ As indicated
in Green's paper, see \cite{mG}, its graded components
play an important r\^ole in getting information about
the minimal free resolution of $M$ over $R.$ Under
additional assumptions on $X$ resp. $M$ there are explicit
geometric interpretations of the Green modules. So e.g. in
the case of $X$ an arithmetically Cohen-Macaulay scheme of
codimension three it follows that
$\mu (I_X) \leq a(I_X) + 1 + \deg X.$

On the other hand the Green modules are intimately related
to the structure of the minimal free resolution of $\mathcal
F_{\bullet}$ of $M.$ In fact they are the
ingredients of a spectral sequence for computing
$\Tor^R_i(K, M), i \in \mathbb Z.$ That is, they describe
in a certain sense the Betti numbers and their shifts.
This is worked out in more detail in Section 5, where
it is shown that
$$
\reg X = \max \{j - i \s j \geq \codim X \mbox{ and } i_j \neq 0 \},
$$
see \ref{5.2}. That is, the regularity of $X$ is completely determined
by the tail of $\mathcal F_{\bullet}.$ More precisely, under
certain additional assumptions there is an explicit
computation of $i_j$ in terms of the graded components
of the Green modules, see \ref{5.5} for the precise statement.
It turns out that this is a generalization of M.~Green's
duality theorem, see \cite[Section 2]{mG}.

On the other hand
Theorem 5.5 is a far reaching generalization of P.~Rao's
observation on how much the resolution of the Hartshorne-Rao
module $M(C)$ of a curve $C \subset \mathbb P^3_K$ determines
the resolution of $\mathcal J_C$ at
the tail, see \cite[(2.5)]{pR}. One application of this type concerns the
resolution of certain curves $C \subset \mathbb P^n_K$ of
arithmetic genus $g_a(C) = 0,$ see \ref{5.7}.

While the applications of our results are motivated by
geometric questions we formulate and prove them in terms
of graded modules over $R$ and its local cohomology modules
$\HH^j(M).$ To this end we fix a few homolo\-gical
preliminaries in Section 2. These concern Koszul
homology, local cohomology, and some basic facts
about spectral sequences. The spectral sequence related
to a double complex is in several variations one of the
basic tools of our investigations. In Section 3 we summarize
the details about the Green modules. The most important
result is \ref{3.4}. It proves the finite length of
$H_i(\lu ; \HH^j(M))$ for a generic system of linear
forms, i.e., for almost all $t \in \mathbb Z$ it follows
that
$$
H_i(\lu ; \HH^j(M))_{i+t} = 0 \, \mbox{ for all } i, j \in \mathbb Z.
$$
Section 4 is
devoted to the estimates of the number of generators of
$I_X,$ i.e., to the desired variations of Dubreil's
theorem. The final Section 5 concerns the relation of
the syzygies of the modules of `deficiencies' $\HH^j(M)$ to those of
$M.$ In particular it yields the new characterization
of $\reg M$ resp. the generalization of Green's duality
theorem.

The author is grateful for the referee's careful reading of the
manuscript.

\section{Koszul homology and local cohomology}
First fix a few notation and conventions. Let $A = \oplus _{n \geq 0}
A_n$ denote a graded Noetherian ring such that $A_0 = K$ is an
infinite field and $A = K[A_1].$ Then $A$ is an epimorphic image
of the polynomial ring $R = K[x_1, \ldots, x_r]$ in the variables
$x_1, \ldots, x_r, \, r = \dim_K A_1.$ Let $M = \oplus_{n \in \mathbb Z}
M_n$ denote a graded $A$-module. For $k \in \mathbb Z$ let $M(k)$
denote the module $M$ with the grading given by $[M(k)]_n = M_{k+n},
\, n \in \mathbb Z.$ Mostly we consider a graded $A$-module as an
module over $R.$ For more details about graded modules and rings see
\cite[1.5]{BH}.

Let $X^{\bullet}, Y^{\bullet}$ denote two complexes of graded $A$-modules.
Let $Z^{\bullet}$ denote the single complex associated to the double
complex $X^{\bullet} \otimes_A Y^{\bullet}.$ Then there are the following
two spectral sequences
\begin{displaymath}
\begin{array}{clllcll}
E_2^{ij} & = & H^i(X^{\bullet} \otimes_A H^j(Y^{\bullet})) & \Rightarrow &
E^{i+j} & = & H^{i+j}(Z^{\bullet}) \quad \mbox{and } \\
`E_2^{ij} & = & H^i(H^j(X^{\bullet}) \otimes_A Y^{\bullet}) & \Rightarrow &
`E^{i+j} & = & H^{i+j}(Z^{\bullet}).
\end{array}
\end{displaymath}
See e.g. \cite[Appendix A3]{dE} or \cite[5.6]{cW} for an introduction and
the basic results concerning spectral sequences.
Here we remark that by the definitions all the homomorphisms are
homogeneous of degree zero.

Let $\lf = f_1, \ldots, f_s$ denote a system of homogeneous elements.
Then $K_{\bullet}(\lf ; A)$ denotes the Koszul complex with  respect
to $\lf.$ Fix the following definitions
\begin{displaymath}
\begin{array}{ll}
K_{\bullet}(\lf ; M) := K_{\bullet}(\lf ; A) \otimes_A M, \quad &
H_i(\lf ; M) := H_i(K_{\bullet}(\lf ; M)), \\
K^{\bullet}(\lf ; M) := \Hom_A(K_{\bullet}(\lf ; A), M), \quad &
H^i(\lf ; M) := H^i(K^{\bullet}(\lf ; M)),
\end{array}
\end{displaymath}
where $i \in \mathbb Z,$ see e.g. \cite[1.6]{BH} or \cite[Section 17]{dE}. Note
that all the modules
resp. complexes are graded. The homomorphims are homogeneous of degree
zero. In particular let $\underline{\mathfrak m} = x_1, \ldots, x_r$ denote a
generating
set
of $\mathfrak m,$ the ideal generated by all forms of positive degree. Then
$K_{\bullet}(\xx ; R)$ provides a finite free resolution of $K,$ the
residue field. Therefore
\begin{displaymath}
H_i( \mathfrak m; M) \simeq \Tor_i^R(K, M) \mbox{ and }
H^i( \mathfrak m; M) \simeq \Ext_R^i(K, M),  i \in \mathbb Z.
\end{displaymath}
In the following split $\mathfrak m$ into two subsets $\xx$ and $\yy.$ Then we
compare their Koszul homologies.

\begin{lemma} \label{2.1} Put $\xx = x_1, \ldots, x_s$ and
$\yy = x_{s+1}, \ldots, x_r$ for an integer $1 \leq s < r.$ Then
$$
\dim _K H_n(\mathfrak m; M) \leq \sum_{i=\max\{0,n-(r-s)\}}^{\min\{s,n\}}
\dim _K H_i(\xx; H_{n-i}(\yy;M)), \quad n \in \mathbb N,
$$
for any finitely generated graded $R$-module $M.$
\end{lemma}

\begin{proof} First note that
$$
K_{\bullet}( \mathfrak m;M) \simeq K_{\bullet}(\xx
; R) \otimes_R K_{\bullet}(\yy; M)
$$
as follows by view of the construction of the Koszul complex. But then there
is the following spectral sequence
$$
E^2_{ij} = H_i(\xx; H_j(\yy; M)) \Rightarrow E_{i+j} =
H_{i+j}(\mathfrak m ;M).
$$
Moreover note that all the $E^2_{ij}$-terms are finite dimensional
$K$-vector spaces. This follows because all of them are annihilated by
$\mathfrak m.$
The subsequent terms $E^n_{ij}$ are subquotients of $E^2_{ij}.$ So they are
also finite dimensional and
$$
\dim_K E^n_{ij}\, \leq \, \dim_K E^2_{ij} \quad \mbox{for all }  n \geq 2.
$$
Now for large $n$ one has $E^{\infty}_{ij} = E^n_{ij}.$ Furthermore
it is known that
$$
E_{i+j} = H_{i+j}(\mathfrak m;M)
$$
admits a finite filtration
whose quotients are $E^{\infty}_{i,n-i},\, i = 0,1, \ldots, n.$ This implies
that
$$
\dim_K H_n(\mathfrak m;M) = \sum_{i=0}^n \dim_K E^{\infty}_{i,n-i}.
$$
So the claim
follows because of the above estimates.
\end{proof}

For further investigations the case of $n = 1$ is of a particular interest.
To this end formulate it as a separate Corollary.

\begin{corollary} \label{2.2} Let $\xx, \yy,$ and $M$ as in \ref{2.1}.
Then
$$
\dim_K H_1(\mathfrak m; M) \leq \dim_K H_0(\xx ;H_1(\yy; M)) +
\dim_K H_1(\xx; H_0(\yy; M)).
$$
\end{corollary}

The general idea behind \ref{2.1} is a bound of the Betti numbers
$$
b_n(M) := \dim _K \Tor_n^R(K, M).
$$
For $n = 1$ and $M = R/I, I$ a homogeneous
ideal of $R,$ this yields a bound for the minimal number of
generators $\mu(I)$ of $I,$ see \ref{4.1}.

For several investigations we need the local cohomology modules $\HH^i(M),
i \in \mathbb Z, $ of $M$ with respect
to $\mathfrak m.$  To this end denote by $K^{\bullet}_f(A)$
the complex $0 \to A \to A_f \to 0,$ where $f$ denotes
a homogeneous element and $A_f$ is the
localization with respect to $f.$ The middle homomorphism denotes
the canonical map into the localization. For
$\underline{\mathfrak m} = x_1. \ldots, x_r$ define
$$
K^{\bullet} := \otimes_{i=1}^r K^{\bullet}_{x_i} \quad \mbox{and } \quad
K^{\bullet}(M) := K^{\bullet} \otimes_A M
$$
the \v Cech complex of $A$ and $M.$ Then there are canonical isomorphisms
$$
\HH^i(M) \simeq H^i(K^{\bullet} \otimes_A M)\, \mbox{ for all }
i \in \mathbb Z.
$$
For the details of this facts see e.g. \cite{rH} or \cite[3.5]{BH}.

\section{The Green modules}
Let $M$ denote a finitely generated graded $R$-module. In this section we
introduce certain invariants related to the local cohomology and the Koszul
homology of $M.$ To this end we need the notion of a generic system of
elements.

\begin{definition} \label{3.1} A system of linear elements
$\lu = l_1, \ldots,l_s$ is said to be a generic linear
system of elements with respect to $M$ provided
$$
l_i \not\in \mathfrak p \mbox{ for all } \mathfrak p \in (\Ass_R(M/(l_1,
\ldots, l_{i-1})M) \setminus \{\mathfrak m\}.
$$
Here $\mathfrak m$ denotes the ideal generated by the variables $x_1, \ldots,
x_r$ in $R.$
\end{definition}

Note that $\lu$ is a generic system of linear elements
if and only if the following quotients
$$
((l_1,\ldots,l_{i-1})M :_M l_i)/(l_1,\ldots,l_{i-1})M, \, i = 1,\ldots,s,
$$
are graded $R$-modules of finite length. This observation is helpful in order
to check whether a given $\lu$ is a generic linear system.
The most important property of a general linear system is related to a certain
finiteness property of $H_i(\lu ; \HH^j(M))$ which we shall prove in this
section. In order to do that we need another auxiliary statement.

\begin{lemma} \label{3.2} Let $\lu = l_1, \ldots, l_s$ denote a generic linear
system with respect to $M.$ Then $H^i(\lu ; M)$ is an $R$-module of finite
length in the following two cases:
\begin{itemize}
\item[(a)] $i < s,$
\item[(b)] for all $i \in \mathbb Z,$ provided $s \geq \dim_R M.$
\end{itemize}
\end{lemma}

\begin{proof} First prove the claim in (b). To this end note that
$$
(\lu, \Ann_R M) H^i(\lu; M) = 0 \, \mbox{ for all }\, i \in \mathbb Z.
$$
Therefore $\Supp _R H^i(\lu; M)
\subseteq \Supp_R M/\lu M.$ Since $s\geq \dim_R M$ it follows
that $\Supp _R M\lu M
\subseteq V(\mathfrak m).$ Recall that $\lu$ is
generically chosen. This proves (b) since $H^i(\lu; M)$
is a finitely generated $R$-module.

The statement in (a) will be shown by an induction on
$d := \dim_RM.$ First note
that the case $d = 0$ is covered by the claim proved in (b).
So let $d > 0.$
First suppose that $\depth_R M > 0.$ Then $l = l_1$ is an $M$-regular
element. The short exact sequence
$$
0 \to M(-1) \to M \to M/lM \to 0
$$
induces short exact sequences
$$
0 \to H^i(\lu;M) \to H^i(\lu; M/lM) \to H^{i+1}(\lu; M)(-1) \to 0
$$
for all $i \in \mathbb Z.$ Note that $l H^i(\lu; M) = 0.$
Hence the induced maps on the Koszul cohomology are trivial.
Now put $\lu' = l_2, \ldots, l_s.$ Then
$$
H^i(\lu;M/lM) \simeq H^i(\lu';M/lM) \oplus H^{i-1}(\lu';M/lM)(-1),
\quad i \in \mathbb Z,
$$
see \cite[1.6]{BH}.
Note that $l$ acts trivially by multiplication on $M/lM,$
Now by induction hypothesis $H^i(\lu'; M/lM)$ is of
finite length for all $i < s-1.$ Therefore
$H^i(\lu; M/lM)$ is of finite length for all $i < s-1.$ Whence the
above short exact sequence proves the claim.

Finally let $\depth_R M = 0.$ Then
$N := \cup_{n \geq 1} (0 :_M l^n)$ is an $R$-module of
finite length as follows by the definition of $\lu.$
Then $\depth_R M/N > 0. $ By the first part of the
inductive step and $\dim _R M/N = d$ the claim is
true for $M/N.$ Note that $\lu$ forms a generic
system of linear forms with respect to $M/N$ as
easily seen by a localization argument with respect to non-maximal
prime ideals. By (b) the
claim is true for $N$ and all $i \geq 0$ since
$\dim_R N = 0.$ So the final statement for $M$
follows from the induced long exact Koszul
cohomology sequence derived from
$0 \to N \to M \to M/N \to 0.$
\end{proof}

In the following consider the functor
$\Hom_K(\Box, K) = {\Box}^{\vee}$ on the category
of graded $R$-modules. By the graded version of the Local Duality Theorem,
see \cite[3.6.19]{BH}, it turns out that there is a natural graded isomorphism
of degree zero
$$
\HH^i(M) \simeq (\Ext_R^{r-i}(M, R(-r))^{\vee}, \quad i \in \mathbb Z.
$$
Put $K^i_M := \Ext_R^{r-i}(M, R(-r)), \, i \in \mathbb Z.$ Then $K^i_M = 0$ for
$i<0$ and $i > \dim_R M := d.$ In particular $K_M := K_M^d$ is called the
canonical module
of $M.$ These modules are studied in a systematic way in \cite[\S 3]{pS2}. Here
we mention only that
$$
\dim_R K^i_M \leq i,\, \mbox{ for } 0 \leq i < d,\, \mbox{ and }\, \dim_R K_M =
d,
$$
see \cite[3.1.1]{pS2} for the details.

For the next results we need another definition of genericity. It is related to
the modules of `deficiency' $K^i_M.$

\begin{definition} \label{3.3} A generic system of linear
elements $\lu = l_1, \ldots, l_s$ is
called a strongly generic linear system of elements
with respect to $M$
provided it is a generic linear system of elements
for all $K^i_M, \, i = 0, \ldots,d.$
\end{definition}

Because $K$ is an infinite field it is clear that
strongly generic linear systems of elements with
respect to $M$ always exist. Their construction is
just an application of prime avoidance arguments.

\begin{theorem} \label{3.4} Suppose that
$\lu = l_1, \ldots,l_s$ is a strongly generic
linear system of elements with respect to $M.$
Then $H_i(\lu; \HH^j(M))$ is a graded $R$-module
of finite length in the following two cases:
\begin{itemize}
\item[(a)] $i < s,$ and
\item[(b)] for all $i \in \mathbb Z$ provided $s \geq j.$
\end{itemize}
\end{theorem}

\begin{proof} First observe that there are canonical isomorphisms
$$
H_i(\lu; N^{\vee}) \simeq (H^i(\lu;N))^{\vee}, \quad i \in \mathbb Z.
$$
This follows because of the isomorphism of complexes
$$
K_{\bullet}(\lu;R) \otimes_R (N)^{\vee} \simeq (\Hom_R (K_{\bullet}(\lu;R),
N))^{\vee},
$$
which is well known. Here $N$ denotes an arbitrary graded $R$-module. Put
$N = K^i_M.$ Then it follows that
$$
H_i(\lu; \HH^j(M)) \simeq (H^i(\lu; K^j_M))^{\vee} \, \mbox{ for all } \, i, j
\in \mathbb Z.
$$
By \ref{3.2} it is known that $H^i(\lu; K^j_M)$ is an $R$-module
of finite length for $i < s$ resp. for all $i \in \mathbb Z$
provided $s \geq \dim_R K^j_M.$ Because of
$\dim_R K^j_M \leq j$ this finishes the proof.
\end{proof}

In his paper \cite{mG} M.~Green considered the following situation. Let
$\mathcal F$ denote a coherent sheaf on $X,$ a compact complex manifold. Then
he considered the vector spaces $\mathcal K^i_{p,q}(X, \mathcal F).$
Let $i \geq 1.$ Then it is easy to see, see \cite{mG}, that
$$
\mathcal K^i_{p,q}(X, \mathcal F) \simeq H_p(\mathfrak m; \HH^{i+1}(M))_{p+q},
$$
where $M$ denotes the associated graded module to $\mathcal F.$ So in an
obvious way the Koszul homology modules
$ H_p(\mathfrak m; \HH^{i+1}(M))$ are graded
analogues of the invariants introduced by M.~Green. As an application of
\ref{3.4} it turns out that $\mathcal K^i_{p,q}(X, \mathcal F) = 0$ for all $q
\ll 0$ resp. for all $q \gg 0$ provided $\lu$ is strongly generically chosen.
For the numerical influence of these finitely many nonvanishing $\mathcal
K^i_{p,q}(X, \mathcal F)$ on free resolutions see the results in Section 5.

The most important feature of $H_p(\mathfrak m; \HH^{i+1}(M))$ is that it is
one of the ingredients of a spectral sequence. In the following let $M$ denote
a finitely generated graded $R$-module. Choose $\lu = l_1,\ldots,l_s, \, s \geq
\dim_R M,$ a generic linear system of elements with respect to $M.$ Then
consider the following complexes $K^{\bullet},$ the \v Cech complex,
$K_{\bullet}(\lu; M),$ the Koszul complex of $M$ with respect to $\lu,$ and
$C^{\bullet} := K^{\bullet} \otimes_R K_{\bullet}(\lu; M).$ Then there is the
following spectral sequence
$$
\HH^i(H_j(\lu ; M)) \Rightarrow H_{j-i}(C^{\bullet}).
$$
Because of the choice of $\lu$ it turns out that $H_j(\lu;M) \simeq
H^{s-j}(\lu;M)(-s)$ are $R$-modules of finite length for all $j \in \mathbb Z,$
see \ref{3.2}. Because of the basic properties of local cohomology it yields
that
$$
\HH^i(H_j(\lu;M)) =
  \begin{cases} H_j(\lu;M) & \mbox{ for } i = 0 \\
                   0       & \mbox{ for } i \not= 0.
  \end{cases}
$$
Therefore the spectral sequence degenerates partially to the isomorphisms
$H_j(C^{\bullet}) \simeq H_j(\lu ; M)$ for all $j \in \mathbb Z.$ The second
spectral sequence for the corresponding double complex is
$H_j(\lu; \HH^i(M)) \Rightarrow H_{j-i}(C^{\bullet}).$ Putting this together it
proves the following

\begin{lemma} \label{3.5} Let $M$ be a finitely generated graded $R$-module.
Let $\lu = l_1,\ldots,l_s, \, s \geq \dim_R M$ be a generic linear system with
respect to $M.$ Then there is a spectral sequence
$$
E_2^{-j,i} = H_j(\lu; \HH^i(M)) \Rightarrow E^{-j+i} = H_{j-i}(\lu; M),
$$
where all the derived homomorphisms are homogeneous of degree zero.
\end{lemma}

In the more special situation of a strongly generic linear system with respect
to $M$ not only $H_{j-i}(\lu;M)$ but also $H_j(\lu; \HH^i(M))$ are modules of
finite length for all $i, j \in \mathbb Z,$ see \ref{3.4}. Therefore there is
an estimate for the length of $H_{j-i}(\lu;M).$

\begin{corollary} \label{3.6} Suppose that
$\lu = l_1,\ldots,l_s, \, s \geq \dim_R M,$ denotes a
strongly generic linear system of elements with respect
to $M.$ Then
$$
L_A(H_n(\lu;M)) \leq \sum_{i=0}^{\min\{s-n,d\}} L_A(H_{n+i}(\lu;\HH^i(M)),
$$ for all $n \in \mathbb Z,$ where $d = \dim_RM.$
\end{corollary}

\begin{proof} First note that $\HH^i(M) = 0$ for
all $i > d$ resp. $H_j(\lu;M) = 0$ for all $j > s.$
Then the estimate follows by the same line of
reasoning as in the proof of \ref{2.1}.
\end{proof}

The spectral sequence in \ref{3.5} has several more
applications in Section 4
and Section 5. Here we want to add just two
simple consequences. They are helpful also in different situations.

\begin{corollary} \label{3.7} Let $M$ be a finitely generated graded
$R$-module. Suppose that $\lu = l_1,\ldots,l_s, \, s \geq \dim_R M$ denotes a
generic linear sytem with respect to $M.$ Then
\begin{itemize}
\item[(a)] $H_{s-t}(\lu; M) \simeq H_s(\lu; \HH^t(M)),$ \, where $t =
\depth_RM,$ and
\item[(b)] $H_{i-d}(\lu;M) \simeq H_i(\lu; \HH^d(M)),$ \, for all $i \in
\mathbb Z,$ provided $M$ is a $d$-dimensional Cohen-Macaulay module.
\end{itemize}
\end{corollary}

\begin{proof} In order to prove (a) consider the spectral sequence in
\ref{3.5}. Take the terms $E_2^{-j,i}$ with $j-i = s-t.$  Then
$$
E_2^{-j,i} =
  \begin{cases}
      0 & \mbox{ for }\, j > s \mbox{ or } i < t \text{ and } \\
      H_s(\lu; \HH^t(M)) & \mbox{ for }\, j = s \mbox{ and } i = t.
  \end{cases}
$$
But this means that the spectral sequence degenerates partially to the desired
isomorphism.

The claim in (b) follows by a similar argument since
$\HH^i(M) = 0$ for all $i \not= d$ in the case of a
Cohen-Macaulay module $M.$
\end{proof}

For an extension of the results of this section to the
situation of a finitely generated module over a local ring, see \cite{pS3}.

\section{Bounds on the number of generators}
For homogeneous ideals $I \subset R = K[x_1,\ldots,x_r], r \geq 2,$ such that
$I$ is a perfect ideal of codimension two it is known that
$$
\mu (I) \leq a(I) + 1,
$$
where $a(I) = \min \{ n \in Z \s I_n \not= 0\},$
the initial degree of $I.$ Note that $a(I)$ is equal
to the minimal degree of a non-zero form contained
in $I.$ This estimate is a generalization of a
corresponding bound given by P.~Dubreil in
the case of $r = 3,$ see \cite{pD}. For the proof see e.g.
\cite{DGM} resp. \cite{pS1}. An approach related to
Hilbert functions is developed in \cite{DGM}, while \cite{pS1}
contains a proof based on the Hilbert-Burch Theorem.

In the following put $\xx = x_1, x_2,\, \yy = x_3,\ldots,x_r, \, r \geq 3,$
where $x_1,\ldots,x_r$ denotes a set of generators of $\mathfrak m.$ In a
certain
sense the following result is a generalization of Dubreil's Theorem.

\begin{theorem} \label{4.1} Let $I \subset R$ denote a
homogeneous ideal of codimension at least two. Then
$$
\mu (I) \leq a(I) + 1 + \mu(H_1(\yy; R/I)),
$$
where $\yy$ is chosen generically with respect to $R/I.$
\end{theorem}

\begin{proof} First put $S = R/\yy R$ and  $J = I S.$
Then we obtain the bound
$$
\mu(I) \leq \mu(J) + \dim_K H_0(\xx; H_1(\yy; R/I)),
$$
as follows by \ref{2.2}. By the generic choice of $\yy$
it is known that $a(I) = a(J).$ Now $J$ is a perfect
ideal of codimension two in $S.$ Therefore
$\mu (J) \leq a(J) + 1.$ Finally the dimension of
the vector space $H_0(\xx; H_1(\yy; R/I))$ coincides
with the number of generators of $H_1(\yy; R/I).$
\end{proof}

In fact \ref{4.1} is a generalization of
J.~Migliore's result, see \cite[Corollary 3.3]{jM},
in the case of the defining ideal $\mathcal J_C$ of a
curve $C \subset \mathbb P^3_K.$ Here we extend his
result to an arbitrary projective scheme.

\begin{corollary} \label{4.2} Let $I \subset R$
denote a homogeneous ideal with
$\codim I \geq 2.$ Put $t = \depth R/I.$ Then
$$
\mu(I) \leq a(I) + 1 + \mu(H_{t+1}(\yy; \HH^t(R/I))),
$$
where $\yy$ is chosen strongly generic with respect to $R/I.$
\end{corollary}

\begin{proof} By \ref{3.7} there is the following isomorphism
$$
H_1(\yy; R/I) \simeq H_{t+1}(\yy; \HH^t(R/I)).
$$
Therefore the claim is a consequence of \ref{4.1}.
\end{proof}

In the situation of $I$ the saturated defining ideal
of curve $C \subset \mathbb P^3_K$ it follows that
$t = 1.$ Therefore $H_2(l_1,l_2; \HH^1(R/I))$ is just
the submodule of
$H^1_{\ast}(\mathcal J_C)$ annihilated by $l_1, l_2,$ see
\cite[Corollary 3.3]{jM}.

Besides of its vanishing it is known that Koszul homology
is difficult to handle. So for the rest of this section
there are several approaches in order
to estimate the term $\mu(H_1(\yy; R/I))$ in \ref{4.1}.

\begin{corollary} \label{4.3} Let $I \subset R$ denote
a homogeneous ideal of codimension at least two with
$d = \dim R/I.$ Then
$$
\mu(I) \leq a(I) + 1 + \sum_{i=0}^d L_R(H_{i+1}(\yy; \HH^i(R/I))).
$$
Moreover, suppose that
$\HH^i(R/I)$ are graded $R$-modules of finite length
for $ i = 0,1,\ldots,d-1.$ Then
$$
\mu(I) \leq a(I) + 1 + \sum_{i=0}^{d-1} \binom{r-2}{i+1}
L_R(\HH^i(R/I)) + L_R(H_{d+1}(\yy; \HH^d(R/I))).
$$
Here $\yy $ is chosen strongly generic with respect to $R/I.$
\end{corollary}

\begin{proof} Under the additional assumption that $\yy $
is a strongly generic system of linear forms with respect
to $R/I$ it follows that $H_{i+1}(\yy; \HH^i(R/I))$ are
graded $R$-modules of finite length, see \ref{3.4}. Then
the spectral sequence in \ref{3.5} provides the following estimate
$$
L_R(H_1(\yy; R/I)) \leq \sum_{i=0}^d L_R(H_{i+1}(\yy; \HH^i(R/I))).
$$
By virtue of \ref{4.1} this proves the first part of the claim.

Under the additional assumption of the finite length of
$\HH^i(R/I)$ for $ i = 0,1,\ldots,d-1,$ it is easy to see that
$$
L_R(H_{i+1}(\yy; \HH^i(R/I))) \leq \binom{r-2}{i+1}
L_R(\HH^i(R/I)), \quad i = 0,\ldots,d-1.
$$
To this end consider the definition of the Koszul homology.
Therefore the second bound follows.
\end{proof}

Note that \ref{4.3} was shown by H.~Martin and J.~Migliore, see
\cite[Theorem 2.5]{MM}, under the additional assumption
that $\Proj R/I$ is equidimensional and a Cohen-Macaulay scheme.
Of particular interest is the case of $\codim I = 2.$
In this situation the term $H_{d+1}(\yy; \HH^d(R/I))$
does not occur since
$d + 1 = r - 1 > r -2,$ the number of elements of $\yy.$

In the following let $\sigma (N)$ denote the socle
dimension of $N.$ That means
$\sigma (N) = \dim_K \Hom_R(R/\mathfrak m, N)$ for
an arbitrary $R$-module $N.$

\begin{corollary} \label{4.4} Let $I \subset R$ denote
a perfect homogeneous ideal of codimension at
least three. Then
$$
\mu(I) \leq a(I) + 1 + \sigma (H^{d+1}(\yy; K_{R/I})),
$$
where $\yy$ is chosen strongly generic with respect to $R/I.$
Here $K_{R/I}$ denotes the canonical module of $R/I.$
\end{corollary}

\begin{proof} Because $R/I$ is a Cohen-Macaulay ring we have to estimate
$\mu(H_{d+1}(\yy; \HH^d(R/I), d = \dim R/I,$ see \ref{4.2}. But now
$$
H_{d+1}(\yy; \HH^d(R/I)) \simeq (H^{d+1}(\yy; K_{R/I}))^{\vee}.
$$
Therefore the dimension of $R/\mathfrak m \otimes_R H_{d+1}(\yy; \HH^d(R/I))$
is equal to the socle dimension of $H^{d+1}(\yy; K_{R/I})$ as easily seen.
\end{proof}

Of a particular interest is the case of a Gorenstein
ideal of codimension three. In this situation it follows:

\begin{corollary} \label{4.5} Let $I \subset R, \yy$ be
as in \ref{4.4}. Suppose that $R/I$ is a Gorenstein ring and
$\codim I = 3.$ Then $\mu(I) \leq 2 a(I) + 1.$
\end{corollary}

\begin{proof} Because $R/I$ is a Gorenstein ring it is known that
$R/I \simeq K_{R/I}.$ Put $\zz = x_3,\ldots,x_{r-1}, \, y = x_r.$
Now define $S = R/\zz R, \, J = IS.$ Then $H^{d+1}(\yy; K_{R/I}) \simeq
S/(J, yS).$ But now the socle dimension of $S/(J, yS)$ is equal to the type
of $S/(J, yS),$ or what is the same, to the minimal number of generators of $L$
minus one, $\mu(L) - 1,$ where $T = S/yS$ and $L = J T.$ Recall that $L$
is a perfect ideal of codimenssion two in $T.$ But then $\mu (L) \leq a(L)
+1$ by Dubreil's Theorem. Finally $a(I) = a(L)$ since $\yy$ is
chosen generically. Therefore by \ref{4.4} the claim is shown to
be true.
\end{proof}

Note that this result follows also by the Buchsbaum-Eisenbud structure
theorem for Gorenstein ideals of codimension three. For the details
see \cite{pS1}. A further result including the degree is the following:

\begin{corollary} \label{4.6} Let $I \subset R, \, \yy$ be as in \ref{4.4}.
Suppose $\codim I = 3.$ Then
$$
\mu(I) \leq a(I) + 1 + e(R/I),
$$
where $e(R/I)$ denotes the multiplicity of $R/I.$
\end{corollary}

\begin{proof} First note that by 4.4 it is obviously true that
$$
\sigma(H^{d+1}(\yy; K_{R/I})) \leq L_R(H^{d+1}(\yy; K_{R/I}))
\leq L_R(K_{R/I}/\zz K_{R/I}).
$$
Here let $\yy$ be generated by $y_1,\ldots,y_{d+1}$ and $\zz = y_1,\ldots,y_d.$
Then $\zz$ forms a system of parameters
for $K_{R/I}$ and $R/I$ as well. Furthermore  $L_R(K_{R/I}/\yy K_{R/I}) \leq
L_R(K_{R/I}/\zz
K_{R/I}).$ Because $R/I$ is a Cohen-Macaulay ring, $K_{R/I}$ is a
Cohen-Macaulay module and therefore
$$
L_R(K_{R/I}/\zz K_{R/I}) = e(\zz; K_{R/I}) = e(\zz; R/I).
$$
Because of the generic choice of the linear elements in
$\zz$ this completes the proof.
\end{proof}

The bound in \ref{4.6} is rather rough. It would be of some interest to
find a common generalization of \ref{4.5} and \ref{4.6}.

\section{Koszul homology and syzygies}
As before let $R = K[x_1,\ldots,x_r]$ denote the polynomial ring in $r$
variables. For a graded $R$-module $M$ define
$$
a(M) = \min \{n \in \mathbb Z \s M_n \not= 0 \} \text{ and }
e(M) = \max \{n \in \mathbb Z \s M_n \not= 0 \}.
$$
It is well known that $e(\HH^i(M)) < \infty$ for all $i \in \mathbb Z.$

\begin{definition} \label{5.1} The Castelnuovo-Mumford regularity
$\reg M$ of $M$ is defined by
$$
\reg M = \max \{ e(\HH^i(M)) + i \s i \in \mathbb Z\}.
$$
Note that $e(0) = - \infty.$
\end{definition}

It is a well known fact that
$$
\reg M = \max\{ e(\Tor_i^R(K,  M)) -i \s  0 \leq i \leq r\}.
$$
So $\reg M$ yields a bound on the maximal degree in a minimal generating
set of the syzygy modules of $M.$ It reflects the structure of the
minimal free resolution $F_{\bullet}$ of $M$ over $R,$ where
$$
F_{\bullet} : 0 \to F_s \to \ldots \to F_i \to \ldots \to F_0 \to M \to 0,
$$
with $F_i \simeq \oplus_{j \in \mathbb Z} R^{i_j}(-j)$ and $i_j = \dim_K
\Tor_i^R(K, M)_j.$ Suppose that $M$ is a Cohen-Macaulay module. Then
$\reg M = e(\Tor_c^R(K, M)) -c,$ where $c = r - \dim_R M$ denotes
the codimension of $M.$ This follows easily since $\Hom_R(F_{\bullet}, R(-r))$
gives a minimal free resolution of  $K_M = \Ext_R^c(M, R(-r)),$ the
canonical module of $M.$

On the other hand it was observed by P.~Rao, see \cite[(2.5)]{pR},
that in the case of $I$ the defining ideal of a curve $C \subset
\mathbb P^3_K$ the Hartshorne-Rao module $M(C) \simeq \HH^1(R/I)$
gives certain information on the tail of the minimal free
resolution of $R/I.$

In the following we shall generalize both of these observations. Firstly
we describe $\reg M$ in terms of the $\Tor$s in a certain range. Secondly
we shall clarify how the minimal free resolutions of $\HH^i(M),$ the
`modules of deficiency', determine the minimal free resolution of $M.$ Both
considerations turn out by a careful study of the spectral sequence
given in \ref{3.5}.

\begin{theorem} \label{5.2} Let $M$ denote a finitely generated
graded $R$-module. Let $s \in \mathbb N.$ then the following two
integers coincide
\begin{itemize}
\item[(a)] $\max \{e(\HH^i(M)) + i \s 0  \leq i \leq s \}$ and
\item[(b)] $\max \{e(\Tor_j^R(K, M)) - j \s r - s \leq j \leq r\}. $
\end{itemize}
In particular for $s = \dim_R M$ it follows that
$$
\reg M = \max \{e(\Tor_j^R(K, M)) - j \s c \leq j \leq r\},
$$
where $c = r - \dim_R M$ denotes the codimension of $M.$
\end{theorem}

Before proving 5.2 we separate two partial results as
Lemmas. They concern results in this direction which seem to be of
some independent interest.

\begin{lemma} \label{5.3} Suppose that $H_s(\mathfrak m; M)_{s+t} \not=
0$ for a certain $t \in \mathbb Z$ and $r - i \leq s \leq r.$ Then there
exists an $j \in \mathbb Z$ such that $0 \leq j \leq i$ and
$\HH^j(M)_{t-j} \not= 0.$
\end{lemma}

\begin{proof} Assume the contrary, i.e., $\HH^j(M)_{t-j} = 0$ for
all $0 \leq j \leq i.$ Then consider the spectral sequence
$$
[E_2^{-s-j,j}]_{t+s} = H_{s+j}(\mathfrak m; \HH^i(M))_{t+s}
\Rightarrow [E^{-s}]_{t+s} = H_s(\mathfrak m;M)_{t+s}
$$
as defined in \ref{4.5}. Recall that all the homomorphisms are
homogeneous of degree zero. Now the corresponding $E_2$-term is a
subquotient of
$$
[\oplus \HH^j(M)^{\binom {r}{s+j}}(-s-j)]_{t+s}.
$$
Let $j \leq i.$ Then this vectorspace is zero by the assumption
about the local cohomology. Let $j > i.$ Then $s + j > s + i \geq r$
and $\binom{r}{s+j} = 0.$ Therefore the corresponding $E_2$-term
$[E_2^{-s-j,j}]_{t+s}$ is zero for all $j \in \mathbb Z.$ But then also
all the subsequent stages are zero, i.e., $[E_{\infty}^{-s-j,j}]_{t+s} = 0$
for all $j \in \mathbb Z.$ Therefore $[E^{-s}]_{t+s} = H_s(\mathfrak m;M)_{t+s}
= 0,$ contradicting the assumption.
\end{proof}

The second partial result shows that a certain non-vanishing of $\HH^i(M))$
yields the existence of a minimal generator of a higher syzygy
module.

\begin{lemma} \label{5.4} Suppose that there are integers $s, b$ such that
the following conditions are satisfied:
\begin{itemize}
\item[(a)] $\HH^i(M)_{b+1-i} = 0$ for all $i < s$ and
\item[(b)] $H_r(\mathfrak m; \HH^s(M))_{b+r-s} \not= 0$
\end{itemize}
Then it follows that $H_{r-s}(\mathfrak m;M)_{b+r-s} \not= 0.$
\end{lemma}

Note that the condition (b) in \ref{5.4} means that
$\HH^s(M)$
possesses a socle generator in degree $b - s.$ Recall that $r$ denotes the
number of generators of $\mathfrak m.$

\begin{proof} As above we consider the spectral sequence
$$
E_2^{-r,s} = H_r(\mathfrak m; \HH^s(M))
\Rightarrow E^{-r+s} = H_{r-s}(\mathfrak m;M)
$$
in degree $b+r-s.$
The subsequent stages of $[E_2^{-r,s}]_{b+r-s}$ are derived by the
cohomology of the following sequence
$$
[E_n^{-r-n,s+n-1}]_{b+r-s} \to [E_n^{-r,s}]_{b+r-s} \to
[E_n^{-r+n,s-n+1}]_{b+r-s}
$$
for $n \geq 2.$ But now $[E_n^{-r-n,s+n-1}]_{b+r-s}$ resp.
$[E_n^{-r+n,s-n+1}]_{b+r-s}$ are subquotients of
$$
H_{r+n}(\mathfrak m; \HH^{s+n-1}(M))_{b+r-s} = 0 \text{ resp. }
H_{r-n}(\mathfrak m; \HH^{s-n+1}(M))_{b+r-s} = 0.
$$
For the second module recall that it is a subquotient of
$$
[\oplus \HH^{s-n+1}(M)^{\binom{r}{r-n}}(-r+n)]_{b+r-s} = 0,
\quad n \geq 2.
$$
Therefore $[E_2^{-r,s}]_{b+r-s} = [E_{\infty}^{-r,s}]_{b+r-s} \not= 0$
and
$$
[E^{-r+s}]_{b+r-s} \simeq H_{r-s}(\mathfrak m; M)_{b+r-s} \not= 0
$$
as follows by the filtration with the corresponding $E_{\infty}$-terms.
\end{proof}

\begin{proof} {\it (Theorem 5.2).} First of all let us introduce two
abbreviations.
Put $a := \max \{ e(\Tor_j^R(K, M)) -j \s r-s \leq j \leq r\}.$ Then by
\ref{5.3} it follows that $a \leq b,$ where $b:= \max \{ e(\HH^i(M)) + i
\s 0 \leq i \leq s\}.$ On the other hand choose $j$ an integer $0 \leq	j
\leq s$ such that $b = e(\HH^j(M)) + j.$ Then $\HH^j(M)_{b-j} \not= 0,$
$ \HH^j(M)_{c-j} = 0$ for all $c > b,$ and $\HH^i(M)_{b+1-i} = 0 $ for all
$i < j.$ Recall that this means that $\HH^j(M)$ has a socle generator in
degree $b - s.$ Therefore Lemma \ref{5.4} applies and
$\Tor^R_{r-j}(K, M)_{b+r-j} \not= 0.$ In other words, $b \leq a,$ as
required.
\end{proof}

An easy byproduct of our investigations is the above mentioned fact
that
$$
\reg M = e(\Tor_c^R(K, M)) - c, \, c = r - \dim M,
$$
provided $M$ is a Cohen-Macaulay module.

\begin{theorem} \label{5.5} Let $M$ be a finitely generated graded $R$-module
with $d = \dim _R M.$
Suppose there is an integer $j \in \mathbb Z$ such that for all $q \in
\mathbb Z$ either
\begin{itemize}
\item[(a)] $\HH^q(M)_{j-q} = 0$ or
\item[(b)] $\HH^p(M)_{j+1-q} = 0$ for all $p < q $ and $\HH^p(M)_{j-1-q} = 0$
for all $p >q.$
\end{itemize}
Then for $s \in \mathbb Z$ it follows that
\begin{itemize}
\item[(1)] $\Tor_s^R(K, M)_{s+j} \simeq \oplus_{i=0}^{r-s} \Tor_{s+i}^R(K,
\HH^i(M))_{s+j}$ provided $s > c,$ and
\item[(2)] $\Tor_s^R(K, M)_{s+j} \simeq \oplus_{i=0}^{d-1} \Tor_{s+i}^R(K,
\HH^i(M))_{s+j} \oplus \Tor_{c-s}^R(K, K_M)^{\vee}_{r-s-j},$
provided $s \leq c,$
\end{itemize}
where $K_M = \Ext_R^c(M, R(-r)), \, c = \codim M,$ denotes the canonical
mo\-dule of $M.$
\end{theorem}

\begin{proof} As above consider the spectral sequence
$$
E_2^{-s-i,i} = H_{s+i}(\mathfrak m; \HH^i(M))
\Rightarrow E^{-s} = H_{s}(\mathfrak m;M)
$$
in degree $s+j,$ see \ref{3.5}. Firstly we claim that $[E_2^{-s-i,i}]_{s+j}
\simeq [E_{\infty}^{-s-i,i}]_{s+j}$ for all $s \in \mathbb Z.$
Because $[E_2^{-s-i,i}]_{s+j}$ is a subquotient of
$$
[\oplus \HH^i(M)^{\binom{r}{s+i}}(-s-i)]_{s+j}
$$
The claim is true provided $\HH^i(M)_{j-i} = 0.$ Suppose that
$\HH^i(M)_{j-i} \not= 0.$
In order to prove the claim in this case too
note that $[E_{n+1}^{-s-i,i}]_{s+j}$ is the cohomology at
$$
[E_n^{-s-i-n,i+n-1}]_{s+j} \to [E_n^{-s-i,i}]_{s+j} \to
[E_2^{-s-i+n,i-n+1}]_{s+j}.
$$
Then the module at the left resp. the right is a subquotient of
$$
H_{s+i+n}(\mathfrak m; \HH^{i+n-1}(M))_{s+j} \text{ resp. }
H_{s+i-n}(\mathfrak m; \HH^{i-n+1}(M))_{s+j}.
$$
Therefore both of them vanish. But this means that the $E_2$-term coincides
with the corresponding $E_{\infty}$-term. So the target of the
spectral sequence $H_s(\mathfrak m;M)_{s+j}$ admitts a finite filtration
whose quotients are $H_{s+i}(\mathfrak m; \HH^i(M))_{s+j}.$ Because
all of these modules are finite dimensional vectorspaces it follows that
$$
H_s(\mathfrak m;M)_{s+j} \simeq \oplus_{i=0}^{r-s} H_{s+i}(\mathfrak m;
\HH^i(M))_{s+j}
$$
for all $s \in \mathbb Z.$

In the case of $s > c$ it is known that $r - s < d.$ Hence the first part
of the claim is shown to be true. In the remaining case $s \leq c$ the
summation is taken from $i = 0,\ldots,d.$ Therefore we have to interprete
the summand $H_{s+d}(\mathfrak m; \HH^d(M))_{s+j}.$ By the Local Duality
Theorem $\HH^d(M) \simeq (K_M)^{\vee}.$ Therefore there are the
following isomorphisms
$$
H_{s+d}(\mathfrak m; (K_M)^{\vee})_{s+j} \simeq
(H^{s+d}(\mathfrak m; K_M)^{\vee})_{s+j} \simeq
H_{r-d-s}(\mathfrak m; K_M)^{\vee}_{r-s-j},
$$
which proves the second part of the claim.
\end{proof}

As an application of \ref{5.5} we derive M.~Green's duality theorem
\cite[Section 2]{mG}, see also \cite[Theorem 1.2]{NP} for a similar approach of
the original
statement.

\begin{corollary} \label{5.6} Suppose there exists an integer $j \in \mathbb Z$
such that
$$\HH^q(M)_{j-q} = \HH^q(M)_{j+1-q} = 0
$$
for all $q < \dim_R M.$ Then
$$
\Tor_s^R(K, M)_{s+j} \simeq \Tor_{c-s}^R(K, K_M)^{\vee}_{r-s-j},
$$
for all $s \in \mathbb Z,$ where $c = \codim M.$
\end{corollary}

\begin{proof} It follows that the assumptions of Theorem 5.5 are satisfied
for $j$ because of $\HH^p(M)_{j-1-p} = 0$ for all $p > \dim M.$ Therefore the
isomorphism is a consequence of (1) and (2) in 5.5. To this end recall that
$$
\Tor^R_{s+i}(K, \HH^i(M))_{s+j} \simeq H_{s+i}(\mathfrak m; \HH^i(M))_{s+j} =
0, $$
as follows by the vanishing of $\HH^i(M)_{s+j}$ for all $j \in \mathbb Z.$
\end{proof}

M.~Green's duality theorem in \ref{5.6} relates the Betti numbers of $M$
to those of $K_M.$ Because of the strong vanishing assumptions in
\ref{5.6} very often it does not give strong information about Betti
numbers. Often it says just the vanishing which follows also by different
arguments, e.g., the regularity of $M.$

Theorem \ref{5.5} is more subtle. In a certain sense it is an extension of
P.~R.~Rao's argument, see \cite[(2.5)]{pR}. We shall illustrate its usefulness
by the following example.

\begin{example} \label{5.7} Let $C \subset \mathbb P^n_K$ denote a reduced
integral non-degenerate curve over an algebraically closed field $K.$ Suppose
that $C$ is non-singular and of genus $g(C) = 0.$ Let $A = R/I$ denote its
coordinate ring, i.e., $R = K[x_0,\ldots,x_n]$ and $I$ its homogeneous defining
ideal. Then
$$
\Tor_s^R(K, R/I)_{s+j} \simeq \Tor^R_{s+1}(K, \HH^1(R/I))_{s+j}
$$
for all $s \geq 1$ and all $j \geq 3.$ To this end recall that $A$ is a
two-dimensional domain. Moreover it is well-known that $\HH^q(R/I) = 0$ for all
$q \leq 0$ and $q > 2.$ Furthermore it is easy to see that $\HH^1(R/I)_{j-1} =
0$ for all $j \leq 1.$ Moreover $\HH^2(R/I)_{j-1-2} = 0$ for all $j \geq 3$ as
follows because of $g(C) = 0.$ That is for $j \geq 3$ one might apply
\ref{5.5}. In order to conclude we have to show that $\Tor^R_{c-s}(K,
K_{R/I})_{r-s-j} = 0$ for $j \geq 3.$ To this end note that
$$
H_{c-s}(\mathfrak m; K_{R/I})^{\vee}_{r-s-j} \simeq H_{s+2}(\mathfrak m;
\HH^2(R/I))_{s+j}
$$
as is shown in the proof of \ref{5.5}. But this vanishes for $j \geq 2$ as is
easily seen.
\end{example}


\begin{thebibliography}{9999}

\bibitem[1]{BH} {\sc W.~Bruns, J.~Herzog:} `Cohen-Macaulay rings', Cambr.
Univ. Press, 1994.

\bibitem[2]{DGM} {\sc E.~D.~Davis, A.~V.~Geramita, P.~Maroscia:}
{\sl Perfect homogeneous ideals: Dubreil's Theorems revisted}, Bull. Sc.
Math. de France, 2$^e$ Ser., {\bf 108} (1984), 143-185.

\bibitem[3]{pD} {\sc P.~Dubreil:} {\sl Sur quelques propri\'et\'es dans le
plan et des courbes gauches alge\'ebriques}, Bull. Sci. math. France
{\bf 61} (1933), 258-283.

\bibitem[4]{dE} {\sc D.~Eisenbud:} `Commutative Algebra (with a view towards
algebraic geometry)', Springer-Verlag, 1995.

\bibitem[5]{mG} {\sc M.~Green:} {\sl Koszul homology and the geometry of
projective varieties}, J. Diff. Geometry {\bf 19} (1984), 125-171.

\bibitem[6]{rH} {\sc A.~Grothendieck:} `Local cohomology', notes by R.
Hartshorne, Lect. Notes in Math., {\bf 41}, Springer-Verlag, 1967.

\bibitem[7]{MM} {\sc H.~M.~Martin, J.~C.~Migliore:} {\sl Submodules of the
deficiency modules and an extension of Dubreil's theorem}, Peprint, 1995.

\bibitem[8]{jM} {\sc J.~C.~Migliore:} {\sl Submodules of the deficiency
module}, J. London Math. Soc., 2$^{nd}$ Ser., {\bf 48} (1993), 396-414.

\bibitem[9]{dM} {\sc D.~Mumford:} `Lectures on Curves on an Algebraic
Surface', Ann. of Math. Studies No. {\bf 59}, Princeton University
Press, 1966.

\bibitem[10]{NP} {\sc U.~Nagel, Y.~Pitteloud:} {\sl On graded Betti numbers
and geometrical properties of projective varieties}, Manuscripta math.
{\bf 84} (1994), 291-314.

\bibitem[11]{pR} {\sc P.~Rao:} {\sl Liaison among Curves in $\mathbb P^3$},
Invent. math. {\bf 50} (1979), 205-217.

\bibitem[12]{pS1} {\sc P.~Schenzel:} {\sl A note on Dubreil's theorems on the
number of generators of perfect ideals of codimension 2}, C. R. Math.
Rep. Acad. Sci. Canada {\bf 6} (1984), 11-14.

\bibitem[13]{pS2} {\sc P.~Schenzel:} `Dualisierende Komplexe in der
lokalen Algebra und Buchs\-baum-Ringe', Lect. Notes in Math., {\bf 907},
Springer-Verlag, 1982.

\bibitem[14]{pS3} {\sc P.~Schenzel:} {\sl On Koszul homology and
local cohomology of modules over local rings}, in preparation.

\bibitem[15]{cW} {\sc C.~Weibel:} `An Introduction to Homological Algebra',
Cambr. Univ. Press, 1994.

\end{thebibliography}
\end{document}